\documentclass[a4paper,11pt]{amsart}
\usepackage{amssymb,amsfonts,amsmath}
\usepackage{layout}
\usepackage[latin1]{inputenc}
\usepackage{amsthm}
\def\d{\delta}
\def\m{\mathcal}

\def\C{\mathbb{C}}
\def\c2{\mathbb{C}^2}
\def\R{\mathbb{R}}
\def\N{\mathbb{N}}

\def\N{\mathbb{N}}

\def\1{\bold{1}}
\def\B{\mathbb{B}}

\def\f{\varphi}

\def\p{\psi}

\newcommand \W {\Omega}

\newcommand \mE {\mathcal E}
\newcommand \mF {\mathcal F}

\newcommand \vk {\chi}
 \newcommand \Ri{ \Rightarrow }
\newcommand \Sub {\Subset}
\newcommand \sub {\subset}
\newcommand \sm {\setminus}

\newcommand \mrm {\mathrm}

\newcommand \ep {\varepsilon}

\newtheorem{lem}{Lemma}[section]

\newtheorem{pro}[lem]{Proposition}
\theoremstyle{definition}
\newtheorem{defi}[lem]{Definition}
\theoremstyle{plain}
\newtheorem{def/not}[lem]{Definition/Notations}
\numberwithin{equation}{section}
\newtheorem{thm}[lem]{Theorem}
\newtheorem{cor}[lem]{Corollary}

\newenvironment{proof3.1}
{\noindent {\it{Proof of theorem 3.1}}}{$\Box$ \linebreak[4]}

\begin{document}

\title[Weighted pluricomplex energy]
{Weighted pluricomplex energy}
\author{ Slimane BENELKOURCHI }
\address{ Institut de Mathématiques de Toulouse, UMR CNRS 5219,
 Universit\'e Paul Sabatier Toulouse 3,
 118 route de Narbonne,
 F-31062 TOULOUSE Cedex 09 (FRANCE)}
\email{benel@math.ups-tlse.fr} 
\maketitle
\begin{abstract}
We study the complex Monge-Ampère operator on the classes of finite
 pluricomplex energy $\m E_\chi (\W)$ in the general case ($\chi(0)=0$
 i.e. the total Monge-Ampère mass may be infinite). We establish an
 interpretation of these classes in terms of the speed of decrease
 of the capacity of sublevel sets and give a complete description of
the range of the operator $(dd^c \cdot)^n$ on the classes $\m
E\chi(\W).$
\end{abstract}

\section{Introduction}
Let $\W \subset \C^n$ be a bounded hyperconvex domain, i.e. a
connected,
 bounded open such that there exists a negative plurisubharmonic $\rho $
such that $\{z\in \W;\ \rho(z)< -c\}\Sub \W, \ \forall c>0$. Such a
function $\rho $
 is called an
 exhaustion function.
 We let
$PSH(\W)$ denote  the cone of  plurisubharmonic functions (psh for
short) on
 $\W$ and
 $PSH^-(\W)$ denote  the
subclass of negative functions.

In two seminal papers \cite{Ce 1}, \cite{Ce 2}, U.Cegrell has
introduced and studied the complex Monge-Amp\`ere operator $(dd^c
\cdot)^n$ on special classes of unbounded plurisubharmonic functions
in $\Omega$, called {\it Energy classes.}  In   \cite{BGZ 1}, a
formalism developed in \cite{GZ}, was used to give a unified
treatment of all these classes in the case of finite total
Monge-Ampère mass. Here, we continue our study in a more general
context. Given an increasing function    $\chi : \R^- \to \R^-,$ we
consider the set ${\mathcal E}_{\chi}(\Omega)$ of plurisubharmonic
functions of finite $\chi$-weighted Monge-Amp\`ere energy. These are
the functions $u \in PSH(\Omega)$ such that there exists a
decreasing sequence $u_j \in \m E_0(\W)$ with limit $u$ and
$$
\sup_{j \in \N} \int_{\Omega} (-\chi) \circ u_j (dd^c u_j)^n
<+\infty,
$$
where ${\mathcal E}_0(\Omega)$ is  the   cone of all bounded
plurisubharmonic functions $\f$ defined on the domain  $\W$ with
finite total Monge-Ampère mass and  $\lim _{z\to \zeta} \f (z)=0,$
for every $\zeta \in \partial \W$. When $\chi(t) = -(-t)^p$ (resp.
$\chi(t) = -1-(-t)^p$), $\mE_\chi(\W)$ is the class $\mE^p(\W)$
(resp. $\mF^p(\W)$ ) studied by  U.Cegrell in \cite{Ce 1}.

The classes ${\mathcal E}_{\chi}(\Omega)$ have very different
properties,
 depending
on whether $\chi(0)=0$ or $\chi(0) \neq 0$, $\chi(-\infty)=-\infty$
or $\chi(-\infty) \neq -\infty$, $\chi$ is convex or concave. If the
function $\chi $ is convex, or concave, then the class $\m
E_\chi(\W)$ is subset of a natural family of psh functions
introduced by U. Cegrell in
\cite{Ce 3} (cf section 4). In particular, we have \\

\noindent{{\bf Proposition A. } \it
  Let $\chi: \R^- \rightarrow \R^-$ be a convex, or concave, increasing
function such that $\chi(-\infty)=-\infty$ and $\chi(0) = 0$. Then
$$
{\mathcal E}_{\chi}(\Omega) \subset {\mathcal N}^{a}(\Omega).
$$
In particular the Monge-Amp\`ere measure $(dd^c u)^n$ of a function
 $u \in \mE _\chi (\W )$
is well defined and does not charge pluripolar sets. More precisely,
}
$$
\mE_\chi (\W)=\left\{ u \in {\mathcal N}(\Omega) \, / \, \chi \circ
u \in L^1((dd^c u)^n) \right\}.
$$
Many properties follow from an interpretation of these classes in
terms of speed of decrease
of the capacity of sublevel sets:\\

\noindent{{\bf Proposition B. } \it If $\chi $ is an increasing
convex function, then we have  }
$$
{\mathcal E}_{{\chi}}(\W)=\left\{ \f \in PSH^-(\W) \, / \,
\int_0^{+\infty} t^n \chi'(-t) \mrm {Cap}_\W(\{\f<-t\}) dt<+\infty
\right\}.
$$
 Here $\mrm {Cap}_\W(\cdot)$
 denotes the Monge-Ampère capacity introduced by E. Bedford
and B.A. Taylor \cite{BT 1}. This yields in particular several
properties: the classes $\m E _\chi(\W)$ are convex, stable under
taking the maximum.

In section 5, we study the range of the complex Monge-Ampère
operator on the classes $\m E_\chi (\W)$ in the case when the
function $\chi $ is convex. Given a positive Borel measure $\mu $ on
$\W,$ we have: \\

\noindent{\bf Theorem C.} {\it Let $\chi : \ \R^- \to \R^- $ be an
increasing convex  function such that $\chi(-\infty)
 = -\infty .$ Then
 there exists  a unique function $\f \in \m E_\chi(\W)$ such
that $\mu = (dd^c \f )^n$ if and only if  there exists a constant
$C>0$ such that
$$
\int_\W -\chi \circ u d\mu \le C_2  \max \left (1,
  \left(\int_\W -\chi
\circ u\left (dd^c u \right )^n\right)^{\frac{1}{n}}\right ) ,
 \ \forall \ u\in \m E _0(\W).
 $$
}

The proof of this theorem remains valid when $\chi(t) = -(-t)^p$ for
$ p>0,$ which yields a simple proof of the main theorem in
 \cite{Ce 1}.

In section 6, using results from \cite{CKZ} and \cite{Ph}, we prove
that, for  almost all weights $\chi $, the functions of the classes
$\m E _\chi (\W)$ admit global subextension with logarithmic growth
and local subextension with finite $\chi$-energy.

\section*{ Acknowledgments} I would like to express my hearty thanks to Ahmed
Zeriahi for many inspiring discussions and for his help in the proof
of Theorem 3.4.   Next my thanks go to Per {\r{A}}hag and the
referees for their valuable suggestions and remarks.
\section{The  class $ \mF(\W)$}
In this section we give some properties of the U.Cegrell class $
\mF(\W)$. The main tool will be the capacity estimate of the
sublevel sets of psh functions. The Monge-Ampère capacity has been
introduced and studied by E.Bedford and A.Taylor in \cite{BT 1}.
Given $K\sub \W$ a Borel subset, its Monge-Ampère capacity
relatively to $\W$ is defined by
$$
\mrm{Cap}_\W (K) := \sup \left \{ \int _K(dd^c u)^n ; \ u \in
PSH(\W),
 \ -1\le u \le 0 \right \}.
$$
Let recall some U.Cegrell's classes.
 The class $\m E (\W)$
 is the set of plurisubharmonic functions
 $u $
such that for all $z_0 \in \Omega $,
 there exists a neighborhood $V_{z_0}$ of
$z_0$ and $u_j \in {\mathcal E}_0(\Omega)$ a decreasing sequence
which converges towards $u$ in $V_{z_0}$ and satisfies $\sup_j
\int_{\Omega} (dd^c u_j)^n <+\infty$. U.Cegrell has shown \cite{Ce
2} that the operator $(dd^c \cdot )^n$ is well defined on
$\mE(\Omega)$ and continuous under decreasing limits. The class
$\mE(\Omega)$ is stable under taking maximum and it is the largest
class with these properties (Theorem 4.5 in \cite{Ce 2}). The class
$\mE(\Omega)$ has been further characterized by Z.Blocki \cite{Bl
1}, \cite{Bl 2}.

The class ${\mathcal F}(\Omega)$ is the ``global version'' of
 $\mE(\Omega)$:
a function $u$ belongs to ${\mathcal F}(\Omega)$ iff there exists
 a decreasing sequence
 $u_j \in {\mathcal E}_0(\Omega)$
converging towards $u$ {\it in all of } $\Omega$, which satisfies
$\sup_j \int_{\Omega} (dd^c u_j)^n<+\infty$.

The class ${\mathcal F}^a(\Omega)$ is the set of functions
 $u \in {\mathcal F}(\Omega)$
whose Monge-Amp\`ere measure $(dd^c u)^n$ is absolutely continuous
with respect to capacity i.e. it does not charge pluripolar sets.
Similarly, ${\mathcal E}^a(\Omega)$ is the set of functions $u \in
\m E(\W)$ whose Monge-Ampère measure
 $(dd^c u)^n$
vanishes on pluripolar sets.

\begin{lem}\label{est}
Fix $\f \in {\mathcal F}(\W)$. Then  for all $s>0$ and $ t > 0$,

\begin{equation}\label{cap}
 t^n
Cap_\W(\f\ <-s-t ) \leq \int_{(\f< -s)} (dd^c \f)^n \leq s^n
Cap_\W(\f < -s ).
\end{equation}
 Therefore $u\in \mF(\W)$ if and only if
$\limsup_{s\to 0} s^n Cap_\W (u < -s)<+\infty.$ In particular, if
 $u\in \mF(\W)$ then
$$
 \int_\W(dd^c u )^n = \lim _{s\to 0} s^n Cap_\W (u < -s)
 $$
and
$$
 \int_{(u=-\infty)}(dd^c u )^n = \lim _{s\to +\infty } s^n
Cap_\W (u <
 -s).
 $$
\end{lem}
Note that the complex Monge-Ampère measure of a psh function $u$ on
$\W$ does not charge pluripolar sets if and only if it puts no mass
on $(u=-\infty)$ (cf. \cite{BGZ 2}). Thus $u\in \mF^a(\W) $ if and
only if $\lim _{s\to +\infty } s^n Cap_\W (u < -s)=0.$

 The  right hand inequality in (\ref{cap}) has proved by S.Kolodziej \cite{K 4}
 when $\f $ is  bounded (see also \cite{B} and \cite{EGZ} for the compact setting).
 For the  convenience of the reader we give here a simple proof which uses the
 same idea.
 \begin{proof} Fix  $s, \ t>0$.
 Let $K \subset \{ \f  < -s -t \} $ be a  compact subset. Then
\begin{multline*}
Cap_\W (K) = \int _\W \left ( dd ^c u_K^* \right )^n =
 \int _{ \{ \f  < -s -t \}}
 \left( dd ^c u_K^* \right )^n \\
= \int _{\{ \f  < -s + tu_K^*  \}}
 \left( dd ^c u_K^* \right )^n = \frac{1}{t^n}\int _
{ \{ \f  < v \}}
 \left( dd ^c v  \right )^n ,
\end{multline*}
where $u_K^*$  is the relative extremal function of the compact $K$
and $ v:= -s + tu_K^* .$ It follows from \cite{BGZ 2} that
\begin{multline*}
\frac{1}{t^n} \int _{ \{ \f  < v \}}
 \left( dd ^c v \right )^n = \frac{1}{t^n}\int _
{ \{ \f  < \max( \f , v) \}}
 \left( dd ^c \max ( \f , v)  \right )^n
  \le\\
 \frac{1}{t^n}\int _{ \{ \f  < \max ( \f,  v) \}}
 \left( dd ^c  \f  \right )^n =    \frac{1}{t^n}\int _{
 \{ \f  <  -s  + t u_K  \}}
 \left( dd ^c  \f)  \right )^n  \le  \frac{1}{t^n}\int _{
 \{ \f  <  -s   \}}
 \left( dd ^c  \f)  \right )^n .
\end{multline*}
Taking the supremum over all K's yields the first inequality. For
the right hand inequality, we have
\begin{multline*}\int _{
 \{ \f  \le   -s   \}}
 \left( dd ^c  \f)  \right )^n = \int _{\W}
 \left( dd ^c  \f)  \right )^n - \int _{\f>-s }
 \left( dd ^c  \f)  \right )^n\\
 = \int _{\W}
 \left( dd ^c \max( \f , -s )  \right )^n - \int _{\f>-s }
 \left( dd ^c  \max( \f , -s )  \right )^n\\ = \int _{\f \le -s }
 \left( dd ^c  \max( \f , -s )  \right )^n \le s^n Cap_\W \{\f \le -s \}.
\end{multline*}

\end{proof}

 It is known (see \cite{Ce 1}, \cite{Ce 2}) that the class $\mF(\W)$
has many properties.
  Namely it
 is a convex cone, stable under  maximum: if $u\in \mF(\W)$
 and $v \in PSH^-(\W)$ then $\max (u, v) \in \mF(\W)$ and  if $u\in \mF(\W)$
 then $\limsup_{u\to \partial \W} u(z)=0$. The subclass $\mF^a (\W)$
 satisfies the same properties.
 All these properties can be deduced easily from
 Lemma \ref{est} using just some basic properties of the Monge-Ampère
 capacity.

 The following corollary generalizes some result in \cite{Dem}.

\begin{cor}
Fix  $u\in \mF(\W )  ,$ and Let $h : ]-\infty , 0] \to ]-\infty , 0]
$ be an increasing function such that $h(0)=0$ and $h\circ u$ is
psh. Then $h\circ u \in \mF(\W )$ if and only if  $h^\prime(0^-) <
\infty.$ Furthermore $ h\circ u \in \mF^a (\W) $ if and only if
$u\in \mF^a(\W) $ or $ h^\prime(-\infty)=0$. Moreover we have
$$
\int_{\W}(dd^ch\circ u )^n =(h^\prime(0^-))^n \int_{\W}(dd^c u)^n.
$$
$$
\int_{(h\circ u=-\infty )}(dd^ch\circ u )^n = (h^\prime(-\infty))^n
 \int_{(u= -\infty)}(dd^c u)^n.
$$
Here $h^\prime(0^-) = \lim_{s\to 0^-} h(s)/s$ and
 $h^\prime(-\infty)=\lim_{s\to +\infty} \frac{h(-s)}{s}.$
\end{cor}
U.Cegrell observed in \cite{Ce 3} that if $u\in \mF(\W) $ then
$-(-u)^{1/n} \not\in \mF(\W) .$ The corollary above  state that
$-(-u)^\alpha \not\in \mF(\W) ,$$\forall \alpha <1.$

We end up this section by extending some result in \cite{Wi}.
\begin{cor}
Let $\W_1$ and $\W_2$ be two hyperconvex domains in $\C^n$ and
$\C^p$
 respectively. Suppose $u_1 \in \mF(\W_1)$ and  $u_2 \in \mF(\W_2),$
then $\max (u_1 , u_2) \in \mF(\W_1 \times \W_2)$ and
$$\int_{\W_1\times \W_2} (dd^c \max (u_1 , u_2))^{n+p} = \int_{\W_1}
 (dd^c u_1)^n  \int_{\W_2}
 (dd^c u_2)^p,
$$
$$\int_{(u_1=-\infty)\times (u_2=-\infty) } (dd^c \max (u_1 , u_2))^{n+p}
 = \int_{(u_1=-\infty)}
 (dd^c u_1)^n  \int_{(u_2=-\infty)}
 (dd^c u_2)^p.
$$

Moreover, $(dd^c \max (u_1 , u_2))^{n+p}$ vanishes on the pluripolar
subsets of $\W_1\times \W_2$ if and only if $ (dd^c u_1)^n $ (or $
(dd^c u_2)^n $ )
  vanishes on the pluripolar subsets of
$\W_1$ (resp. of $\W_2$).
\end{cor}
\begin{proof} Observe that
\begin{multline*}
\{  (z , w) \in \W_1 \times \W_2 ; \max (u_1(z) , u_2(w))
 \le -s \} = \\
  \{  z  \in \W_1 ; u_1(z)  \le -s \}  \times \{  w  \in \W_2 ;
 u_2(w)  \le -s \}.
 \end{multline*}
Then it follows from \cite{Bl 3} that
\begin{multline*}
s^{n + p}\mrm{cap }_{\W_1\times \W_2}(\{  (z , w) \in \W_1 \times
\W_2 ; \max (u_1(z) , u_2(w)
 \le -s \}) =\\
  s^n\mrm{cap }_{\W_1}( \{  z  \in \W_1 ; u_1(z)  \le -s \}) s^p\mrm{cap }
_{\W-2} ( \{w  \in \W_2 ;
 u_2(w)  \le -s \}.
 \end{multline*}
 Hence the desired results follow by Lemma \ref{est}.
\end{proof}
\section{Capacity of sublevel set}
It's well known that if $u \in PSH^-(\W)$ is any psh function then
 for every compact $K\Sub \W $  there is a
constant $C>0$ such that
$$
\mrm{Cap}_\W  ( { \{u <-s\}}\cap K) \le \frac{C}{s }  , \quad
\forall  \
  s>0.
 $$
But if $u \in  \mE(\W)$, the capacity of sublevel set
 decreases at least like
$s^{-n}$, i.e. for every compact $K\Sub \W $  there is a constant
$C>0$ such that
$$
\mrm{Cap}_\W  ( { \{u <-s\}}\cap K) \le \frac{C}{s^n }  , \quad
\forall  \
  s>0.
 $$
In fact this is a necessary condition (cf Lemma \ref{est}) but not
sufficient to get
 $u\in \mE(\W).$ Indeed, let $\B\sub \C^n ,\ n\ge 2,$ the unit ball,
  we consider the psh function $u(z)= -(-\log |z_1|)^\frac{1}{n}.$ It's
 clear that it satisfies the last condition but $u \notin \mE(\W),$
 cf \cite{Ce 2}, \cite{Ce 4}.

 In this section, we show that if the capacity of sublevel
set of a psh function $u$ decreases fast enough then its complex
Monge-Ampère $(dd^c u)^n$ is well defined.

Denote by
 $ \m P_n(\W)$ the space
of all negative psh function $u\in PSH^-(\W)$ such that
$$
 \int _0^{\infty} s^{n-1} \mrm{Cap }_\W(\{ u < -s \}\cap K ) ds
 < \infty,
 $$
 for every compact $ K\Sub \W .$

Bedford has introduced the following class (see \cite{Be}). Let
$\theta \ : \R \to \R $ be a decreasing function such that $t\to
-(-t\theta(-t))^{1/n}$ is an increasing and convex function on
$]-\infty , 0]$ and
\begin{equation}\label{bed}
\int_1^{+\infty}\frac{\theta (t)}{t} dt <+\infty.
\end{equation}
Define $\m B(\W)$  to be the class of negative function $u\in
PSH^-(\W)$ such that for any $z_0\in \W$ there exist a neighborhood
$\omega$ of $z_0$, a negative psh  function $\psi$ and a decreasing
function $\theta $ satisfying  (\ref{bed}) such that
   $-(-\psi \theta(-\psi ))^{1/n} \le u $ on $\omega.$
\begin{pro}
 For any hyperconvex domain $\W \Sub \C^n,$ we have
  $\m B(\W) \sub  \m P_n(\W).$ In particular, for any negative psh function
  $v$ on $\W$ and any $0 <\alpha < 1/n,$ $-(-v)^\alpha \in \m
  P_n(\W)$.
\end{pro}
\begin{proof}
    It follows from the definition of   $\m B(\W)$ that, for any $\omega
 \Sub \W $ and $s>0$
    \begin{equation}\label{sub}
    \{u<-s\}\cap \omega \sub \{-(-\psi \theta(-\psi ))^{1/n} <-s\}\cap
    \omega=  \{-\psi \theta(-\psi ) >s^n\}\cap
    \omega .
    \end{equation}
    Let $ \kappa $  be a function such that $ \kappa^{'} = \theta $
    and  $ \kappa(0) = 0. $ The function $ \kappa $ is concave.
    Hence $$ \kappa (t) \ge t   \theta(t) , \  \forall t>0,$$    which
    together with (\ref{sub}) yield
\begin{multline*}
\int _0^{\infty} s^{n-1} \mrm{Cap }_\W(\{u \le -s \}\cap \omega ) ds
\le \int _0^{\infty} s^{n-1} \mrm{Cap }_\W(\{\kappa(- \psi) \ge s^n
\}\cap \omega ) ds  \\
 \le C_1 +   \int _1^{\infty} s^{n-1} \mrm{Cap
}_\W(\{ \psi \le - \kappa^{-1}(s^n) \}\cap \omega ) ds \le \\
C_1 + C_2\int _1^{\infty} s^{n-1} \frac{1}{\kappa^{-1}(s^n)} ds =
C_1+C_2 \int_1^\infty \frac{\theta(t)}{t}dt
 < \infty,
 \end{multline*}
 which completes the proof.
\end{proof}
More generally, let us consider an increasing function $h :\ \R^-
\to \R^-$. Then we have:
\begin{pro} Suppose that  $h$ satisfies
\begin{equation}\label{formule}
\int^{+\infty} \frac{(-h(-s))^{n-1} h^{'}(-s)}{s} ds <+\infty.
\end{equation}
Then for any psh function $u\in PSH^-(\W)$ such that $h\circ u \in
PSH\-(\W)$ we have $h\circ u \in   \m P_n(\W).$ Moreover, if $h$ is
convex, then $h\circ PSH^-(\W) \sub   \m P_n(\W).$
\end{pro}
 The following lemma (cf \cite{CKZ}) will be useful later on.
\begin{lem}\label{zeriahi}
 For any psh function $u \in \mE(\W)$, we have
\begin{equation}\label{zer}
\int _{B}(dd^{c }u)^n \le (\parallel u\parallel _{B})^n \mrm{Cap}_\W
(B ) ,
\end{equation}
provided that  $ \parallel u\parallel _{B} = \sup _B |u|<\infty$.
\end{lem}
\begin{proof} Denote $M = \sup_B |u|<\infty, $ and fix $\ep >0.$ Since $B\sub
\{u>-M-\ep\},$ it  follows from
 \cite{BGZ 2}
$$
\int _{B}(dd^{c }u)^n = \int _{B}(dd^{c }\max(u, -M-\ep))^n <
(M+\ep)^n
 \mrm{Cap}_\W (B).
$$
Letting $\ep \to 0$ yields the desired estimate.
\end{proof}
 Here we will show that the complex Monge-Ampère operator is well defined in
the space $ \m P_n(\W)$ and puts no mass on pluripolar sets.
\begin{thm}\label{choquet}
For every hyperconvex domain $\W \Sub \C^n$, we have
 $$
  \m P_n(\W) \sub \mE^a(\W).
 $$

Conversely, if $u\in \mE(\W )$ then there exists an increasing
convex function $\chi \ :\ \R^- \to \R^-$ such that
$$\int _0^{\infty} s^{n-1}\chi^{'}(-s) \mrm{Cap }_\W(\{z \in K \ ;\ u(z) \le -s \}\ )
ds  < \infty,$$ for all compact $K\Sub \W$.
\end{thm}

\begin{proof}The last statement is an immediate consequence of Corollary 4.4 in \cite{BGZ 2}.
To prove the first one, fix $u\in \m P_n(\W)$. It follows from
\cite{Ce 2} that there exists a decreasing sequence $u_j\in \m E _0
(\W) $ such that $\lim _j u_j =u.$ Let $ B \Sub \W$ be a ball and
consider, for $j\ge 1$,  the function $\tilde{u}_j$ defined by
 \begin{equation}
   \tilde{u}_j(z):= \sup \{ v(z); \ v\in PSH^-(\W) \ and \ v\le u_j
   \ \ in \ B\}
\ \ z\in \W. \label{red}
\end{equation}
It's clear that $ \tilde{u}_j $
 decreases to  $u_B$ defined by
$$u_B(z) =\sup \{ v(z); \ v\in PSH^-(\W) \ and \ v\le u
   \ \ in \ B\}
\ \forall  z\in \W.$$ So, it's enough to prove that
$$\sup _{j} \int _{\W }(dd^c  \tilde{u}_j)^n < \infty.$$
In fact, this a simple consequence of some  precise estimate of the
Monge-Ampère mass in terms of capacity of sublevel set which can be
stated as follows. There exists a constant $C = C(n)$ depending only
in $n$ such that
\begin{equation}  \label{cln}
  \int_K (dd^c \f)^n  \le C \int_0^{+\infty}s^{n-1}
  \mrm{Cap}_\W ({K\cap \{ \f\le -s\}) }ds,
\end{equation}
for any negative bounded psh function $\f$ and any Borel subset
$K\Sub \W.$

Indeed, it follows from Lemma \ref{zeriahi}
\begin{eqnarray*}
\int_K (dd^c \f)^n &=& \sum_{k=-\infty}^{k=+\infty}\int _{K\cap
\{2^{k-1}\le -\f<2^k\} }
 (dd^c \f)^n \\
&\le & \sum_{k=-\infty}^{k=+\infty} 2^{-kn}
\mrm{Cap}_\W ({K\cap \{2^{k-1}\le -\f<2^k\} })\\
&\le&
C \sum_{k=-\infty}^{k=+\infty } \int _{2^{k-1}}^{2^k} n s^{n-1}
\mrm{Cap}_\W ({K\cap \{ -\f\ge s \} })\\
&\le & 
C  \int _0^{\infty}  s^{n-1}
\mrm{Cap}_\W ({K\cap \{ \f\le -s\}) }ds  \\
&\le & 
C  \int _0^{\infty}  s^{n-1} \mrm{Cap}_\W ({K\cap \{ \f\le -s\}) }ds
.
\end{eqnarray*}
Now, we apply the estimate (\ref{cln}) to $ \tilde{u}_j$, to get
 $$ \int _{\W }(dd^c  \tilde{u}_j)^n =
 \int _{B }(dd^c  \tilde{u}_j)^n\le
 C  \int _0^{\infty}  s^{n-1} \mrm{Cap}_\W ({K\cap \{ u \le -s\}) }ds
 < \infty.$$
 Which prove that $u_B
\in \m F(\W)$ and therefore $u\in \mE(\W).$ Since the Monge-Ampère
capacity $ \mrm{Cap}_\W (\cdot )$ vanishes on pluripolar sets, it
follows that
 $u_B \in \m F^a(\W) $ and then  $u\in \mE^a(\W).$
\end{proof}
\begin{cor}For any hyperconvex domain $\W \Sub \C ^n,$ we have
$\m B(\W)\sub \m E ^a(\W), $ i.e. for any function $u\in \m B(\W),$
 the complex Monge-Ampère measure $(dd^cu )^n $  is well defined
 and puts no mass
on the pluripolar sets.

If $h: \ \R^- \to \R^-$ is an increasing convex function satisfying
the condition (\ref{formule}), then $h\circ PSH^-(\W) \sub   \m
E^a(\W).$

 In particular, for
any $0<\alpha < 1/n, $  the psh function $ -(-u)^\alpha \in
\mE^a(\W)$.
\end{cor}
A similar result of the second statement, with a deferent proof,
 has been obtained recently by Z.Blocki in \cite{Bl 4}.
 The author wishes to thank to the anonymous referee for sending him the recent paper
 \cite{Bl 4}.

 The first statement has been also proved in \cite{Ce 4} and
 \cite{CKZ}.
\section{The weighted energy class }
\begin{defi}
 Let $\chi : \R^-  \to \R^- $ be an
 increasing
 function. We let $ \m E_\chi (\W ) $
 denote the set of all functions $u \in PSH(\Omega)$ for which
 there exists a sequence $u_j \in {\mathcal E}_0(\Omega)$
 decreasing to $u$ in $\Omega$ and satisfying
 $$
\sup_{j \in \N} \int _ {\W}(-\chi) \circ u_j \, (dd^c u_j )^n
<\infty.
$$
\end{defi}
This definition clearly contains the  classes of U.Cegrell:
\begin{itemize}
\item $\m E_\chi (\W )={\mathcal F}(\W)$ if $\chi$ is bounded and $\chi(0)
 \neq 0$;
\item $\m E_\chi (\W )={\mathcal E}^p(\W)$ if $\chi(t)=-(-t)^p$;
\item $\m E_\chi (\W)={\mathcal F}^p(\W)$ if $\chi(t)=-1-(-t)^p$.
\end{itemize}

Let us stress that the classes ${\mathcal E}_{\chi}(\Omega)$ are
very different whether $\chi (0) \neq 0 $ (finite total
Monge-Amp\`ere mass) or  $\chi (0)=0$, $\chi(-\infty) = -\infty$ or
$\chi(-\infty) \not = -\infty, $ and  $\chi $  is convex or concave.

The case  $\chi (0) \neq 0 $ was studied in \cite{BGZ 2}, here we
consider the general case $\chi (0)=0$.

It is useful in practice to understand these classes through the
speed of decreasing of the capacity of sublevel sets.

\begin{defi}
$$
\hat{{\mathcal E}}_{\chi}(\W ) :=\left\{ \f \in PSH^-(\W) \, / \,
\int_0^{+\infty} t^n \chi'(-t) \mrm {Cap}_\W(\{\f<-t\}) dt<+\infty
\right\}.
$$
\end{defi}

The classes ${\mathcal E}_{\chi}(\Omega)$ and $\hat{{\mathcal
E}}_{\chi}(\Omega)$ are closely related:

\begin{pro}
The classes $\hat{{\mathcal E}}_\chi(\W ) $ are convex and  if $\f
\in \hat{{\mathcal E}}_{\chi}(\W ) $ and $  \p \in PSH^-(\W)$, then
$\max(\f, \p )  \in \hat{{\mathcal E}}_{\chi}(\W)$.

One always has $\hat{{\mathcal E}}_{\chi} (\W ) \subset \mE _\chi
(\W ) $, while
$$
{\mathcal E}_{{\chi}}(\W) \subset \hat{{\mathcal E}}_{\hat{
\chi}}(\W), \text{ where }  \hat{\chi}(t) = \chi(t/2).
$$
\end{pro}
\begin{proof}
Cf Proposition 4.2.  in  \cite{BGZ 2}.
\end{proof}
\begin{cor} Let $\chi : \R^-  \to \R^- $ be an
 increasing
 function.
If $u\in \m E_\chi (\W)$ then  $\limsup_{z\to \zeta } u(z) = 0 , \
\forall \zeta \in \partial \W .$
\end{cor}
\begin{proof}
In fact, we prove the following claim  which has its self interest.

If a subset $E\sub \W$ has a ``big contact''  with the boundary
$\partial \W$
 of $\W$, then its Monge-Ampère Capacity is infinite. For instance, if
$E = B\cap \W$, where $B$ is a ball centered at some point in
$\partial \W$.

Indeed, let $K_j$ be an increasing sequence of regular compact
subsets
 such that
 $E = \cup K_j $. The extremal function $u_{K_j} \in \m E _0(\W)$ and decreases
to the extremal function $u_E$. It's clear that $u_E\not \in
\mF(\W)$. Thus
$$
\sup_j  \mrm {Cap}_\W(K_j)=\sup _j \int_\W (dd^c u_{K_j})^n =
+\infty.
$$
Therefore $\mrm {Cap}_\W(E)=+\infty$.

Now, we prove the corollary. Assume that there exists a $\zeta_0 \in
\partial \W,$ such that  $\limsup_{z\to \zeta_0 } u(z) = \d <0.$
This yields that there exists a small ball centered
 at $\zeta_0$ such that
 $B\cap \W \sub \{u< \delta/2\}.$ Then, it follows from the claim
that  $$\mrm {Cap}_\W \{u< -s\}=+\infty, \ \forall s \le -\d /2,$$
which contradicts the fact $u\in {\mathcal E}_{{\chi}}(\W) \subset
 \hat{{\mathcal E}}_{\hat{
\chi}}(\W).$
\end{proof}
\begin{thm}\label{integration} Let $\chi : \R^-  \to \R^- $ be a
 convex, or concave,
 increasing function such that $\chi(-\infty ) = -\infty $
and $\chi(0) = 0.$  Then
 $$
\m E_\chi (\W) \sub \mE^a(\W).
$$
Moreover,
 if  $u \in \mE(\W) $ and $v\in \mE_\vk (\W)$ are such that
$(dd^cv)^n \le (dd^c u)^n ,$ then  $u \le v.$
\end{thm}
\begin{proof}
Fix $u\in \m E_\chi (\W)$, denote $u_j$ a defining sequence such
that
$$\sup _j \int_\W- \chi (u_j)(dd^c u_j )^n <+\infty .$$

1)° If $\chi $ is convex: It clear that
$$
\sup _j \int_\W- \chi (u_1)(dd^c u_j )^n  \le \sup _j \int_\W- \chi
(u_j)(dd^c u_j )^n <+\infty .$$ So it follows from \cite{Ce 3} that
$u\in \m E(\W).$

2)° If  $\chi $ is concave: It follows from the concavity assumption
that
$$
  -\chi (-j )= \chi (0 ) -\chi (-j ) \le j  \chi^{'} (-j ).
$$
Then  for any compact $K \Sub \W$,
\begin{multline*}
\int_0^{+\infty} t^{n -1} \mrm {Cap}_\W(\{\f<-t\}\cap K ) dt\le\\
 C_1
+C_2 \int_0^{+\infty} t^n \chi'(-t) \mrm {Cap}_\W(\{\f<-t\})
dt<+\infty.
\end{multline*}
Therefore, from Theorem \ref{choquet}, we deduce $u\in  \m E(\W).$

Now we prove the second statement. In fact we will adapt the  same
idea as in the proof of Theorem 6.2
 in \cite{Ce 1} for the case $\mE_p(\W).$
  Fix  $\rho \in \mE_0(\W), $
 not identically $0.$ We assume that $-1 \le \rho <0.$

First, we assume that $\chi $ is convex. Then for any  $j\ge 1,$ we
have
$$
(dd^c \max (v, j\rho))^n = \1_{\{v>j\rho\}} (dd^c v )^n +  \1
_{\{v\le j\rho\}} (dd^c \max ( v, j\rho ) )^n,
$$
where $\1 _A$ is the characteristic function for the set $A.$ By
\cite{K 1}
 there exist $g_j \in \mE_0$ such that $ (dd^c g_j)^n = \1 _{\{v\le
j\rho\}} (dd^c \max ( v, j\rho ) )^n.$ Thus $(dd^c (u+g_j ))^n \ge
 (dd^c \max ( v, j\rho ) )^n.$  It follows from the comparison principle
for bounded psh functions  (see for example \cite{BT 1},  \cite{Ce
2}) that
 $u+g_j \le  \max ( v, j\rho ) )^n .$
Hence
$$
 u + \limsup_{j\to \infty} g_j \le v= \lim_{j\to \infty }\max ( v, j\rho ).
$$
Now it's enough to prove that $ \limsup_{j\to \infty} g_j =0$ a.e.
Denote $w_m := (\sup _{k\ge m } g_k)^*$, we prove that $ \int _\W
-\vk (w_m ) (dd^c w_m )^n =0,$ and this implies that $w_m =0,\ a.e.$
 Let $j\ge m.$ By
integration by part, we have
\begin{multline*}
 \int _\W -\vk (w_m ) (dd^c w_m )^n  \le  \int _\W -\vk (m\rho )
 (dd^c w_j )^n \le
  \int _\W -\vk (m\rho ) (dd^c g_j )^n \\
\le  \sup_{z\in \W} \frac{-\vk (m\rho (z) )}{-\vk (j\rho(z) )}\
 \int _\W -\vk (j\rho )  \chi _{\{v\le
j\rho\}} (dd^c \max ( v, j\rho ) )^n \\
\le
  \sup_{z\in \W} \frac{-\vk (m\rho (z) )}{-\vk (j\rho(z) )}
\  \sup _{j\ge m} \int _\W
  -\chi ((\max (v, j\rho)) (dd^c \max ( v, j\rho ) )^n <+\infty.
\end{multline*}
We claim that
$$
\lim_{j\to \infty} \sup_{z\in \W} \frac{-\vk (m\rho(z) )}{-\vk
(j\rho(z) )}=0.
$$
Indeed, for $z\in \W$, put  $s =\rho(z).$ Assume, on the contrary,
that
\begin{equation} \label{contrad}
\limsup_{j\to \infty} \sup_{-1 \le s \le 0} \frac{-\vk (ms )} {-\vk
(js )}>\delta >0.
\end{equation}
Then there exists a sequence $s_j$ converging towards $0$ such that
$\frac{-\vk (ms_j )} {-\vk (js_j )}>\delta >0. $ Since  $ms_j\to 0,$
as $j\to \infty$,
 it follows that  $js_j\to 0,$ as $j\to \infty.$
Since $\chi $ is convex,
 we have
 $$\frac {-\vk (ms_j )}
{-\vk (js_j )} \thicksim \frac{ms_j}{js_j}=\frac{m}{j} \to 0, \
\text{as}\
 j\to +\infty,
 $$
which contradicts (\ref{contrad}). Therefore,the claim is proved.
 Hence
  $ \limsup_{j\to \infty} g_j = 0.\ a.e. $

Now, if $\chi $ is concave. We modify  slightly the above proof.
Indeed, since  $\chi $ is concave, the function $\chi ^{-1}(j\rho )
\in \m E_0(\W)$ for any $j>0.$ Then
 \begin{multline*}
(dd^c \max (v, \chi^{-1}(j\rho )))^n\\ = \1_{\{v>
\chi^{-1}(j\rho)\}} (dd^c v )^n +  \1 _{\{v\le
 \chi^{-1}(j\rho)\}} (dd^c \max ( v,  \chi^{-1}( j\rho )) )^n.
\end{multline*}
We consider the function $g_j \in \mE_0(\W)$ satisfying
$$ (dd^c g_j)^n
= \1 _{\{v\le \chi ^{-1}(j\rho)\}} (dd^c \max ( v, \chi ^{-1}(j\rho)
) )^n.$$
 Then we repeat
 the same arguments as above.
\end{proof}

Note that if $u\in \m E_\chi (\W) $ is such that
 $\int_\W(dd^c u)^n < +\infty$ then $u\in \mF(\W)$. Therefore, by
 Lemma \ref{est},
 the total mass  $\int_\W(dd^c u)^n$ depends only on the behavior of $u$ near
 $\partial \W.$
 Now, if  $\int_\W(dd^c u)^n = +\infty$ then
 $\int_\W(dd^c \max (u, -s))^n = +\infty,\  \forall s>0, $ and since
$\int_{(u=-s)}(dd^c \max (u, -s))^n <  +\infty$ (cf Lemma
\ref{zeriahi}) it follows that $\int_{(u>-s)}(dd^c u)^n = +\infty\
\text{and} \ \int_{(u\le -s)}(dd^c u)^n <  +\infty, \
 \forall s>0. $
\begin{lem} \label{ceg}
If $u \in \m E_\chi (\W)$ then there exists a decreasing sequence
$u_j \in {\m E_0}(\Omega)$ with $\lim u_j = u$ and
$$
\lim_{j\to \infty }\int_\W  (- \chi)\circ  u_j (dd^c u_j )^n =
\int_\W (- \chi ) \circ u (dd^c u )^n <+\infty.
$$
\end{lem}
This result was proved by U.Cegrell (cf \cite{Ce 1}) for the classes
$\mE_p(\W)$. The same proof still valid  in the general context. For
the convenience  of the reader we give here the proof.
\begin{proof}
It follows from \cite{K 1} that  there exists, for each $ j\in \N $,
a function $u_j \in \mathcal{E}_0 (\W ) $
 such that $(dd^c u_j )^n = {\bf 1}_{\{u>j\rho \}} (dd^c u)^n$, where
 $\rho \in \m E _0 (\W)$ any defining function for $\W = \{\rho <
0\}$. Observe  that $ (dd^c u )^n \ge (dd^c u_{j+1} )^n \ge  (dd^c
u_j )^n$. We infer
  from the comparison principle that  $(u_j) $  is a decreasing sequence
and   $\lim_j u_j= u$. The monotone convergence theorem thus yields
\begin{multline*}
\int_\W  (- \chi)\circ  u_j (dd^c u_j )^n \\= \int _\W (- \chi)
\circ u_j {\bf 1}_{\{u>j\rho \}} (dd^c u)^n \to \int_\W (-  \chi )
\circ u (dd^c u )^n <+\infty.
\end{multline*}
\end{proof}
The following capacity  estimates of sublevel sets
 will be useful later on.
\begin{pro}\label{part}
Let $\chi: \R^- \rightarrow \R^-$ be an increasing convex, or
concave, function such that $\chi(-\infty)=-\infty$ and $\chi(0) =
0$. Then
$$
\mrm {Cap}_\W(\{\f<-2s\})  \le \frac{1}{|s^n \chi(-s)|} \int _{(\f
<-s)}-\chi(\f) (dd^c \f )^n,
$$
for any $s>0$ and any function  $\f \in \m E_\chi (\W) .$
\end{pro}
\begin{proof} Follows from Lemma \ref{est} by approximating  $\f$
by $\f_j\in \mE_0(\W)$ given by the   lemma above.
\end{proof}

\begin{pro}Let $\chi: \R^- \rightarrow \R^-$ be an increasing
convex, or concave, function such that $\chi(-\infty)=-\infty$ and
$\chi(0) = 0$.
 Then there exists a constant $C = C(\chi )$ such
that
$$
\mrm {Cap}_\W(\{\f<-s\})  \le \frac{C}{s^n} \int_\W -\chi
(\frac{\f}{s})(dd^c \f )^n , \ \forall s>0, \forall \f \in \m E_\chi
(\W).
$$
\end{pro}
\begin{proof}
First we give the proof in the case $n=2$. Let $K \Sub \{\f<-s\}$ be
a compact subset, $u_K$ denotes its relative extremal function.
Choose $\chi _1 : \ \R^- \to \R^- $ to be  an increasing function
such that $\chi_1 ^{''} = \chi $ and $\chi_1(0)=0$. Then
\begin{equation}  \label{p1}
dd^c \chi _1(\f) =\chi_1 ^{''}(\f) d\f \wedge d^c \f + \chi_1
^{'}(\f)dd^c \f \le \chi_1 ^{'}(\f)dd^c \f,
\end{equation}
and
\begin{equation} \label{p2}
-dd^c \chi _1^{'}(\f) =-\chi_1 ^{'''}(\f) d\f \wedge d^c \f - \chi_1
^{''}(\f)dd^c \f \le -\chi(\f)dd^c \f.
\end{equation}
It follows from \cite{Ce 2} that there exists a decreasing sequence
$\f_j\in \m E _0(\W) \cap C(\bar{\W})$ such that $\f_j \searrow \f.$
Then  integrating by part together with the previous inequalities
yield
\begin{eqnarray*}
  \int_K (dd^c u_K)^n  &\le &  \int_K \frac{-\chi_1(\f /s)}{-\chi_1(-1)}
(dd^c u_K)^n =
  \lim_j\int_K \frac{-\chi_1(\f_j /s)}{-\chi_1(-1)}(dd^c u_K)^n
  \\
   &=&\lim_j \frac{1}{-\chi_1(-1)} \int_\W -u_Kdd^c{-\chi_1(\f_j /s)}
\wedge(dd^c u_K)^{n-1} \\
   &\le &\lim_j \frac{C}{s} \int_\W -u_K {\chi_1^{'}(\f_j /s)} dd^c{\f_j }
\wedge(dd^c
   u_K)^{n-1}\\
   &\le &\lim_j \frac{C}{s} \int_\W  {\chi_1^{'}(\f_j /s)} dd^c{\f_j }
\wedge(dd^c
   u_K)^{n-1}\\
   & \le & \lim_j \frac{C}{s} \int_\W  u_K dd^c {\chi_1^{'}(\f_j /s)}
\wedge dd^c{\f_j }\wedge(dd^c
   u_K)^{n-2}\\
   &\le & \lim_j\frac{C}{s^2} \int_\W -  {\chi_1^{''}(\f_j /s)}
 (dd^c{\f_j })^2\wedge(dd^c
   u_K)^{n-2} \\
   &=&  \frac{C}{s^2} \int_\W -  {\chi_1^{''}(\f /s)}  (dd^c{\f })^2
\wedge(dd^c
   u_K)^{n-2} .
\end{eqnarray*}
 For the general case, we use the same arguments. Indeed, we
 consider an increasing function   $\chi _1 : \ \R^- \to \R^- $
such that $\chi_1 ^{(n)} = \chi $ and $\chi_1(0)=0$. Then, the
repeated application of inequalities (\ref{p1}), (\ref{p2}) and the
integration by part yields the desired estimate.
\end{proof}
Hereafter,  we will see that in fact, the classes $\mE_\chi (\W )$
live in some natural set of psh functions introduced by U.Cegrell in
\cite{Ce 3}. Let us recall its definition. Let $\W_j\Sub \W$ be an
increasing sequence of strictly pseudoconvex domains such that $ \W=
\cup _j \W_j.$ Let $u\in \mE(\W)$ be given and put
$$
 u_{\W _j}:= \sup \left \{\f \in PSH(\W); \ \f \le u \ \text{on}\
 \W\sm \W_j\right \}.
$$
Then the sequence $u_{\W_j} \in \m E(\W)$ is increasing, so
$\tilde{u}:= (\lim_j u_{\W_j})^* \in \m E(\W).$ The definition of
$\tilde{u}$  is independent of the choice of the sequence $\W_j$ and
is maximal i.e. $(dd^c \tilde{u})^n =0.$ $\tilde{u}$ is the smallest
maximal psh function above $u$. Define $\m N(\W):= \{u\in \mE(\W); \
\tilde{u} =0\}.$ In fact, this class is the analogous of potentials
for subharmonic functions. Also, denote $\m N^a(\W)=\m E^a(\W)\cap
\m N(\W).$
\begin{pro}\label{prop}
Let $\chi: \R^- \rightarrow \R^-$ be a convex, or concave,
increasing function such that $\chi(-\infty)=-\infty$ and $\chi(0) =
0$. Then
$$
{\mathcal E}_{\chi}(\Omega) \subset {\mathcal N}^{a}(\Omega).
$$
In particular the Monge-Amp\`ere measure $(dd^c u)^n$ of a function
 $u \in \mE _\chi (\W )$
is well defined and does not charge pluripolar sets. More precisely,
$$
\mE_\chi (\W)=\left\{ u \in {\mathcal N}(\Omega) \, / \, \chi \circ
u \in L^1((dd^c u)^n) \right\}.
$$
\end{pro}
\begin{proof} Fix $u\in \mE_\chi (\W) $ and $u_j \in \m E _0(\W )$ a defining
sequence
 such that
$$\sup _j \int_\W- \chi (u_j)(dd^c u_j )^n <+\infty .$$
It follows from the upper semi-continuity of $u$ that $ - \chi
(u)(dd^c u )^n $ is bounded
 from above
by any cluster point of the bounded sequence $- \chi (u_j)(dd^c u_j
)^n.$
 Therefore
$\int_\W (-\chi )\circ u (dd^c u )^n < +\infty , $  in particular
$(dd^c u )^n $ does not charge the set $\{\chi(u) =
 -\infty \},$ which
 coincides with $\{u= -\infty \}$, since $\chi(-\infty ) = -\infty .$
 It follows from   Theorem
 2.1 in \cite{BGZ 2},
$ (dd^c u)^n$ does not charge pluripolar sets. Now it remains to
prove that $ u \in {\mathcal N}(\Omega) $ ie. the smallest maximal
function above $u$ is null. Let $\tilde{u} $ be a  such function.
Then $u\le \tilde{u}\le 0$, thus $\tilde{u} \in \hat{{\mathcal
E}}_{\hat{ \chi}}(\W) \sub {{\mathcal E}}_{\hat{ \chi}}(\W).$ It
follows from
 Lemma \ref{ceg} that there exists a decreasing sequence $\tilde u_j
 \in {\m E_0}(\Omega)$
 with $\lim \tilde u_j = \tilde u$ and
$$
\lim_{j\to \infty }\int_\W  (- \chi)\circ  \tilde u_j (dd^c \tilde
u_j )^n = \int_\W (- \chi ) \circ \tilde u (dd^c \tilde u )^n
<+\infty.
$$
Hence Lemma \ref{est} implies that $ \int_0^{+\infty} t^n
\chi'(-t/4) \mrm {Cap}_\W(\{\tilde u <-t\}) dt = 0$, this yields
that $\tilde u = 0 $.

 To prove the last assertion, it remains to show the reverse
inclusion  $$ \mE_\chi (\W)\supset\left\{ u \in {\mathcal N}(\Omega)
\, / \, \chi \circ u \in L^1((dd^c u)^n) \right\}.$$ This is an
immediate consequence of Lemma \ref{ceg}.
\end{proof}

Note that, unlike the case $\chi(0)\not = 0$ with the class
$\mF(\W)$ (cf \cite{BGZ 2}), we have
$$
\bigcap_{\substack{ \vk(0) =0 \\
\chi(-\infty )=-\infty} } \mE_\vk (\W )  \varsubsetneq \m N(\W )\cap
L^\infty(\Omega), \quad \mbox{and} \bigcup_{\substack{\chi(0) = 0,\\
\chi(-\infty)=-\infty}} {\mathcal E}_{\chi}(\Omega) \varsubsetneq\m
N^a(\W ) .
$$
One can see \cite{Ce 3} for examples of functions in the class $\m
N^a(\W )\cap L^\infty(\Omega) $ which do not belong to any $\mE_\vk
(\W
).$\\

Let  $\chi : \R^-  \to \R^- $ be an increasing function. We say that
 $\chi $ is {\it admissible }  if and only if $\chi $ is convex or concave
and  if there
 exists a constant $M >0 $ such that
\begin{equation}\label{quasi}
\chi^{'} (-2 s) \le M \chi^{'} (- s), \ \forall s>0.
\end{equation}
Observe that any homogenous function $\chi (t) = -(-t)^p$ $p\ge 1$,
 is admissible.
Anther example of admissible function which is not homogenous
 (cf \cite{GZ}) is
$ \chi(t) = -(-t)^p(\log (-t + e))^\alpha, $ $p\ge 1$ and $\alpha
>0.$

\begin{pro} \label{bene}
If $\chi $ is an increasing admissible function, then we have
$$
{\mathcal E}_{{\chi}}(\W)=\left\{ \f \in PSH^-(\W) \, / \,
\int_0^{+\infty} t^n \chi'(-t) \mrm {Cap}_\W(\{\f<-t\}) dt<+\infty
\right\}.
$$
\end{pro}
\begin{proof}
 Follows easily from Lemma \ref{est} and (\ref{quasi}).
\end{proof}
\begin{thm}\label{continuity}
Let $\chi: \R^- \rightarrow \R^-$ be an admissible  increasing
function such that $\chi(-\infty)=-\infty$ and $\chi(0)  = 0$. Fix
$u\in \mE_\chi(\W)$ and set  $u^j =  \max (u , -j).$ Then for each
Borel subset $B\subset \W ,$
$$
\lim _{j\to \infty}\int_B (dd^cu^j)^n = \int_B (dd^cu)^n ,
$$
 {and}
$$
\int_B \chi (u^j)(dd^cu^j)^n \to \int_B \chi (u)(dd^cu)^n. $$

Furthermore, if $u_j$ is any decreasing sequence in $\mE_\chi(\W)$
converging to $u$, Then
$$
\lim _j \int_\W \chi (u_j)(dd^cu_j)^n = \int_\W \chi (u)(dd^cu)^n.
$$
\end{thm}
The first statement, as we will see in the proof,  still valid for
all
 weight $\chi$.
\begin{proof}Let  $B \subset \W $ be  a
Borel subset. If $\int _B (dd^c u )^n =+\infty $ then for any $j>0$,
$\int _B (dd^c u^j )^n=+\infty  .$ So we assume that $ \int _B (dd^c
u )^n< +\infty.$  It follows from Lemma \ref{zeriahi} and
Proposition \ref{part}
\begin{multline*}
\left |\int _B (dd^c u^j) ^n -\int _B (dd^c u )^n \right | \le
 \int _{ \{u\le -j\}}
 (dd^c u^j) ^n +
\int _{\{u\le -j\}} (dd^c u )^n \\
\le   j^n \mrm {Cap}_\W(\{ u <-j\}) +  \int_{(u<-j)} \frac{-\chi
({u})}{-\chi (-j)}
(dd^c u )^n \\
\le  \frac{2^{n+1}}{-\chi (-j/2)} \int_{(u<-j/2)}-\chi ({u})(dd^c u
)^n ,
 \to 0 ,\  \mbox{as} \ j\to + \infty.
\end{multline*}
The proof that $\chi \circ u^j (dd^c u^j )^n $ converges strongly
towards
 $\chi \circ u (dd^c u )^n $  goes along  similar lines, first  observe that
 from Lemma \ref{zeriahi}, we have
 \begin{eqnarray}
\int_{\{u\le -j \}} -\chi \circ u^j (dd^c u^j )^n  &=& -\chi (-j
)\int_{\{u\le -j \}}
 (dd^c u^j )^n\\
 &\le &  -\chi (-j ) j^n \mrm {Cap}_\W(\{ u <-j\}). \label{n1}
\end{eqnarray}
Since  $ \chi $ is an admissible function, it follows that there
exists a constant $C>1$ such that
$$
 -\chi (-2s ) \le - C \chi (-s ),\ \forall s>0
$$
This yields
\begin{multline}
 \lim_{j\to \infty }-\chi (-j ) j^n \mrm {Cap}_\W(\{ u <-j\})  \le
 \lim_{j\to \infty }-C  \chi (-j/2 ) j^n \mrm {Cap}_\W(\{ u <-j\})\\
  \le
  \lim_{j\to \infty }-2^{n+1}C  \chi (-j ) j^n \mrm {Cap}_\W(\{ u <-2j\})\\
  \le
 \lim_{j\to \infty }2^{n+1}C \int_{\{u\le -j \}}-\chi (u ) (dd^c u
 )^n=0.     \label{n2}
\end{multline}
Then (\ref{n1})and  (\ref{n2}) together with Proposition \ref{part}
imply that
$$\lim_{j\to +\infty }\int_{\{u\le -j \}} -\chi \circ u^j (dd^c u^j)^n
=0.$$
 Hence the proof of the second statement is completed.

 Now, once the first and second assertions are proved, we apply the same
 proof as that
 of Theorem 3.4 in \cite{BGZ 2} to show the last statement.
\end{proof}

We conclude this section with a characterization of bounded function
in the classes $\mE_\chi(\W)$,
  extending Y. Xing's main result in
\cite{Xi}.
\begin{pro}
Let $u\in \m E_\chi (\W)$. Then $u$ is bounded in the domain $\W$ if
and only if there exist constants $A>0$ and $B$ such that for any
real $k< B $ with $Cap_\W (u<k) \not = 0$ we can find an increasing
sequence $k\le k_1 <k_2<\cdots <k_s =B$ with $k_1<k+1$ and
$$
\sum _{j=2}^{s} \left ( \frac{\int _{(u<k_j)}(dd^c
u)^n}{Cap_\W(u<k_{j-1})}\right )^{1/n} < A.
$$
\end{pro}
\begin{proof}
The necessary implication is obvious.
To show the sufficient one, assume on the contrary that $u$ is
unbounded. Then $Cap_\W(u<k) \not =0$ for all $k<0$. It follows from
Lemma \ref{est}
$$
B-1-k\le \sum _{j=2}^sk_j -k_{j-1} \le \sum _{j=2}^s \left (
\frac{\int _{(u<k_j)}(dd^c u)^n}{Cap_\W(u<k_{j-1})}\right )^{1/n} <
A.
$$
Hence $B -1-k \le A$ for all $k<B$, which is impossible. The proof
is complete.
\end{proof}

\section{The range of the complex Monge-Amp\`ere operator}
The image of the complex Monge-Amp\`ere operator acting on the
classes
 $\m E_p(\W),$
 has been extensively studied by U.Cegrell. The main result of his study,
 achieved in
\cite{Ce 1}, is given as follows. Given a positive measure $\mu ,$
then there exists  a unique function $\f \in \m E_p(\W)$ such that
$\mu = (dd^c \f )^n$ if and only if there exists a  constant $C>0$
such that
\begin{equation} \label{ep}
\int_\W (-u)^p d\mu \le C\left ( \int_\W (-u)^p(dd^c u )^n \right
)^\frac{p}{n+p}, \qquad \forall  u\in \m E _0(\W).
\end{equation}
Observe that this necessary and sufficient condition is equivalent
to the following: The operator $u\to \int_\W (-u)^p d\mu $ is
uniformly  bounded on the compact ``pseudo-ball'' $\tilde{\m
E}_p(\W): = \{u\in {\m E}_p(\W); \  \int_\W (-u)^p(dd^c u )^n\le
1\}.$ The following theorem extends U.Cegrell's main result \cite{Ce
1}.

\begin{thm} Let $\chi : \ \R^- \to \R^- $ be
an increasing convex  function such that $\chi(-\infty)
 = -\infty .$
  The following conditions are equivalent:\\

(1)
 there exists  a unique function $\f \in \m E_\chi(\W)$ such
that $\mu = (dd^c \f )^n$;\\

(2) there exists a  constant $C_1>0$ such that
\begin{equation}\label{quant}
\int_\W -\chi \circ u d\mu \le C_1 , \ \forall \ u\in \tilde{\m E
_0}(\W),
\end{equation}

(3) there exists a  constant $C_2>0$ such that

\begin{equation}\label{quant1}
\int_\W -\chi \circ u d\mu \le C_2  \max \left (1,
  \left(\int_\W -\chi
\circ u\left (dd^c u \right )^n\right)^{\frac{1}{n}}\right ) ,
 \ \forall \ u\in \m E _0(\W).
\end{equation}

Here   $\tilde{\m E _0}(\W): = \{u\in {\m E _0} (\W); \  \int_\W
-\chi \circ u(dd^c u )^n\le 1\}.$
\end{thm}

\begin{proof} We prove that $(1)\Rightarrow (2)\Ri (3) \Ri (1).$\\
We start with $(3) \Ri (1).$ It follows from \cite{BGZ 2} (see also
Proposition \ref{bene}) that the class $\m E_\chi(\W )$
characterizes pluripolar sets.
 Then the assumption (\ref{quant1}) on $\mu $ implies in particular
 that it vanishes on pluripolar sets. It follows from \cite{Ce 2}
 that there exists a function $u \in \m E _0(\W)$ and
 $f \in L_{loc}^1 \big((dd^c u )^n\big)$ such that
$\mu = f (dd^c u )^n. $

Consider $\mu_j:=\min (f, j ) (dd^c u )^n$. This is a finite measure
which is bounded from above by the complex Monge-Amp\`ere measure of
a bounded function. It follows therefore from \cite{K 1} that there
exist $\f_j \in \m E _0 (\W)$
 such that
$$
(dd^c \f_j )^n = \min (f, j ) (dd^c u )^n.
$$
The comparison principle shows that $\f_j$ is a decreasing sequence.
Set  $\f =\lim_{j\to \infty } \f_j$. It follows from (\ref{quant1})
that
$$
\int_\W -\chi (\f_j) (dd^c \f_j)^n \le C_2 \max\left(1,
\left(\int_\W -\chi (\f_j) (dd^c \f_j)^n\right) ^{1/n}\right).
$$
Hence
$$
  \sup_j   \int_\W -\chi (\f_j) (dd^c \f_j)^n \le  C_2^{n/n-1}
  <\infty.
$$
So it follows from Proposition  \ref{bene} that
$$
\sup_j \int_0^{+\infty} t^{n }\chi'(-t) \mrm {Cap}_\W(\{\f_j<-t\})
dt<+\infty ,
$$ which implies that
  $$
  \int_0^{+\infty} t^{n }\chi'(-t) \mrm {Cap}_\W(\{\f<-t\}) dt<+\infty
  .
$$
Then $\f \not \equiv -\infty$
  and  therefore                 $\f \in
 \mE_\vk(\Omega) $.

 We conclude now by continuity of the complex
Monge-Amp\`ere operator along decreasing sequences that
 $(dd^c \f )^n = \mu.$ The unicity of $\f $ follows from the
 comparison principle
(Theorem \ref{integration}).

Now, we prove $(2) \Ri (3)$. Let $\p \in \m E _0(\W),$ denote $E_\vk
(\p): = \int_\W -\chi (\p) (dd^c \p)^n.$
  If $\p \in
\tilde{\m E _0}(\W)$,   i.e.  $E_\vk (\p)  \le 1 $  then
$$
 \int_\W -\chi (\p) d\mu \le C_1.
$$
If  $E_\vk (\p)  > 1 .$ The function $\tilde{\p}$ defined by
$$
\tilde{\p}:=\frac{\p}{E_\vk (\p)^{1/n}} \in   \tilde{\m E _0}(\W).
$$
Indeed, from the monotonicity of $ \chi $, we have
$$
\int_\W -\chi ( \frac{\p}{E_\vk (\p)^{1/n}}) (dd^c \frac{\p}{E_\vk
(\p)^{1/n}})^n \le \frac{1}{E_\vk (\p)} \int_\W -\chi (\p) (dd^c
\p)^n = 1.
$$
It follows from (\ref{quant}) and the convexity of $\chi$
\begin{multline*}
\int_\W -\chi (\p) d\mu  \le
 {E_\vk (\p)^{1/n}} \int_\W
-\chi ( \frac{\p}{E_\vk (\p)^{1/n}}) d\mu  \le  C_1. {E_\vk
(\p)^{1/n}}.
\end{multline*}
Hence we get (3) with $C_2 = \max(1, C_1).$

For  the proof of the remaining implication $ (1) \Ri (2),$ we use
the same
 idea as in \cite{GZ}. Let $u\in \tilde{\m E _0}(\W)$
 and $\f \in \mE_\chi(\W)$.
 Observe that
 for any $s>0$, we have
 $$
 (u<-s ) \sub (u<\f -\frac{s}{2}) \cup (\f <-\frac{s}{2}).
 $$
 Hence
 \begin{multline}
  \int_\W -\chi \circ u (dd^c \f )^n = \int_0^\infty -\chi {'}(-s)
   \int_ {(u<-s)} (dd^c \f )^n ds \\
   \le \int_0^\infty \chi {'}(-s)
   \int_ {(u<\f -\frac{s}{2})} (dd^c \f )^n ds +  \int_0^\infty \chi {'}(-s)
   \int_ {(\f< -\frac{s}{2})} (dd^c \f )^n ds \\
   \le 2\int_0^\infty \chi {'}(-2s)
   \int_ {(u<\f -s)} (dd^c \f )^n ds + 2 \int_0^\infty \chi {'}(-2s)
   \int_ {(\f< -s)} (dd^c \f )^n ds .  \label{E1}
\end{multline}
The convexity  of $\chi $  yields  that \begin{equation} \label{E2}
 \chi {'}(-2s)\le
M \chi {'}(-s), \ \forall s>0.
\end{equation}
It follows by the comparison principle that, for all $s>0$
\begin{equation} \label{E3}
\int_ {(u<\f -s)} (dd^c \f )^n \le \int_ {(u<\f -s)} (dd^c u )^n \le
\int_ {(u< -s)} (dd^c u )^n .
\end{equation}
Together (\ref{E1}), (\ref{E2}) and (\ref{E3}) imply that there
exists a constant $C$ independent of $u$ such that
 $\int_\W -\chi \circ u (dd^c \f)^n \le C , \ \forall u \in
\tilde{\m E _0}(\W). $
\end{proof}

Note that if  $\chi $ is homogenous, i.e. $\chi(t) =- (-t)^p $ with
$p>0,$  then  the above theorem still valid,
 but we replace the assertion (3) by the following

$(3^{'})$ there exists a  constant $C_2^{'}>0$ such that

\begin{equation}\label{quant2}
\int_\W -\chi \circ u d\mu \le C_2^{'}
  \max \left(1, \left(\int_\W -\chi
\circ u(dd^c u )^n\right)^{\frac{p}{n+p}}\right ) ,
 \ \forall \ u\in \m E _0(\W),
\end{equation}
which, thanks to the homogeneity,  is equivalent to (\ref{ep}). In
particular, this generalizes the U.Cegrell's main theorem in
\cite{Ce 1}for $p\ge 1$ and in \cite{ACP} for $0<p\le 1.$
\section{Subextension in the class $\mE\chi$}
Here we will show that functions in the classes $\m E_\chi(\W)$
admit subextension. We need to recall the usual Lelong class of psh
functions. Let $\gamma >0$ be a positive real. Then
$$
\m L_\gamma (\C^n   ) : = \left \{ \f \in PSH(\C^n); \ \limsup_{r\to
+\infty}\frac{\max _{||z||=r} \f (z)}{\log r } \le \gamma \right \}.
$$

\begin{pro}Let $\chi : \ \R^- \to \R^- $ be
an increasing  function such that $\chi(-\infty) = -\infty $ and
$$
\int ^{+\infty}\frac{1}{s|\chi(-s)|^{1/n}} ds <+\infty .
$$
Then for any function $\f \in \mE_\chi (\W)$ and any
 $\ep >0$, there exists a
function $ U_\ep \in \m L_\ep (\C^n)$ such that $U_\ep \le \f$ on
$\W.$
\end{pro}
\begin{proof}Define the function $h(s) :=  \mrm {Cap}_\W(\{ u <-s\})$.
 It follows from the proof of Theorem \ref{continuity} that
$$
 \mrm {Cap}_\W(\{ u <-s\})
 \le \frac{2^n}{s^n|\chi (-s/2 )|}
 \int_{\{u\le -s/2 \}}-\chi (u) (dd^c u )^n
.$$ Then
$$
\int ^\infty  h(s)^{1/n}ds \le 2 \left (  \int_{\W}-\chi (u) (dd^c u
)^n\right )^{\frac{1}{n}}
 \int ^\infty
 \frac{1}{s|\chi (-s/2 )|^{1/n}} ds <+\infty .
$$
Hence the assertion follows from Theorem 4.1 in \cite{CKZ}.
\end{proof}
\begin{thm}
Let $\W\sub \tilde{\W} \sub \C^n$ be hyperconvex domains.
 Let $\chi : \ \R^- \to \R^- $ be
an increasing  function such that $\chi(-\infty) = -\infty .$ If
 $u \in \mE_\chi (\W)$, then there exists $\tilde{u} \in \mE_\chi (\tilde{\W})$
such that $\tilde{u} \le u$ on $\W$, $(dd^c \tilde{u})^n \le (dd^c
u)^n $ on $\W$ and $E_\chi(\tilde{u}) \le E_\chi (u).$
\end{thm}
\begin{proof}
With slightly different notations, the proof is identical to that in
the case $\mE_p(\W). $ We refer the reader to  \cite{Ph} for
details.
\end{proof}



\begin{thebibliography}{widestlabel}

\addcontentsline{toc}{chapter}{Bibliography}


\bibitem {ACP} P.{\r{A}}hag \& R.Czyz R. \& H.H.Pham: Concerning the energy
class $\m Ep$ for $0<p<1$, Ann. Polon. Math. 91 (2007), 119-130.



\bibitem {Be} {E.Bedford}: { Survey of pluripotential theory,} Several
 Complex Variables, Mittag-Leffler institute
1987-88, Math. Notes, 38(1993), 48--95, Princeton Uni. Press,
Princeton.



 \bibitem  {BT 1} E.BEDFORD \& B.A.TAYLOR: A new capacity for
 plurisubharmonic
   functions. Acta Math. {\bf 149} (1982), no. 1-2, 1--40.

  \bibitem  {BT 2} E.BEDFORD \& B.A.TAYLOR:    Fine topology, \v Silov
boundary,
   and $(dd\sp c)\sp n$.  J. Funct. Anal.  {\bf 72}  (1987),  no. 2, 225--251.

 \bibitem  {B}  S.BENELKOURCHI:   A note on the approximation of
 plurisubharmonic
  functions.
C. R. Math. Acad. Sci. Paris, {\bf 342} (2006), 647--650.

\bibitem  {BGZ 1} S.BENELKOURCHI \& V.GUEDJ \& A.ZERIAHI: A priori
 estimates
 for weak solutions of
  complex Monge-Amp\`ere equations., Ann. Scuola Norm. Sup. Pisa
Cl. Sci. (5), Vol. VII(2008), 81--96.


 \bibitem  {BGZ 2} S.BENELKOURCHI \& V.GUEDJ \& A.ZERIAHI:
 Plurisubharmonic functions
 with weak singularities,
Proceedings from the Kiselmanfest, 2006. Acta Universitatis
Upsaliensis, Proceedings of the conference in honor of C.Kiselman
(Kiselmanfest, Uppsala, May 2006) (in press).

 \bibitem  {Bl 1} Z.BLOCKI: On the definition of the
 Monge-Amp\`ere operator in
 $\Bbb C\sp 2$.
  Math. Ann.  {\bf 328}  (2004),  no. 3, 415--423.

 \bibitem  {Bl 2} Z.BLOCKI: The domain of definition of the complex
  Monge-Amp\`ere operator.
   Amer. J. Math.  {\bf 128}  (2006),  no. 2, 519--530.

\bibitem  {Bl 3} Z.BLOCKI: Equilibrium measure of a product subset
 of $\C\sp n$.  Proc. Amer. Math. Soc.  128  (2000),  no. 12, 3595--3599.

  \bibitem  {Bl 4} Z.BLOCKI: Remark on the definition of the complex
  Monge-Ampère operator. Preprint
  http://www.mitag-leffler.se/preprints/0708s/info.php?id=02

 \bibitem  {Ce 1} U.CEGRELL: Pluricomplex energy. Acta Math. {\bf 180}
 (1998), no. 2,
   187--217.

 \bibitem  {Ce 2} U.CEGRELL: The general definition of the complex
 Monge-Amp\`ere operator.
   Ann. Inst. Fourier (Grenoble)  {\bf 54}  (2004),  no. 1, 159--179.

\bibitem  {Ce 3} U.CEGRELL: A general Dirichlet problem for  of the
complex Monge-Amp\`ere operator, Ann. Polon. Math. 94 (2008),
131-147.

 \bibitem  {Ce 4} U.CEGRELL:  Explicit calculation of a Monge-Amp`ere measure,
  Actes des rencontres
d'analyse complexe, 25-28 Mars 1999. Edited by Gilles Raby and
Frédéric Symesak. Atlantique. Université de Poitiers, 2000.

\bibitem  {CKZ}
U.CEGRELL \& S.KOLODZIEJ \& A.ZERIAHI: Subextension of
plurisubharmonic functions with weak singularities.  Math. Z.  {\bf
250}  (2005),  no. 1, 7--22.

 \bibitem  {Dem}J.-P.DEMAILLY: Monge-Amp\`ere operators, Lelong numbers
 and intersection
   theory. Complex analysis and geometry, 115--193, Univ. Ser. Math., Plenum,
   New York (1993).


 \bibitem  {EGZ} P.EYSSIDIEUX \& V.GUEDJ \& A.ZERIAHI:  Singular K\"ahler-Einstein
  metrics. J. Amer.
Math. Soc., to appear.


 \bibitem {GZ} V.GUEDJ \& A.ZERIAHI: The weighted Monge-Amp\`ere energy
   of quasiplurisubharmonic functions.  J. Funct. Anal.  250  (2007),
 no. 2, 442--482.

 \bibitem {K 1}S.KOLODZIEJ: The range of the complex Monge-Amp\`ere operator.
    Indiana Univ. Math. J. {\bf 43} (1994), no. 4, 1321--1338.

 \bibitem  {K 2} S.KOLODZIEJ: The complex Monge-Amp\`ere equation. Acta Math. {\bf 180}
   (1998), no. 1, 69--117.

 \bibitem  {K 4} S.KOLODZIEJ: The complex Monge-Amp\`ere equation and
 pluripotential
  theory.
   Mem. Amer. Math. Soc.  {\bf 178}  (2005),  no. 840, x+64 pp.

 \bibitem {Ph}  H.H.PHAM: Pluripolar sets and the subextension in
 Cegrell's classes,  Complex Var. Elliptic Equ.  53  (2008),  no. 7, 675--684.



\bibitem{Wi} J.WIKLUND: Topics in pluripotential theory.
Doctoral Thesis No. 30, 2004, Umea university, Sweden.



 \bibitem  {Xi} Y.XING: Complex Monge-Ampère measures of
 plurisubharmonic functions with bounded values near the boundary.
  Canad. J. Math.  52  (2000),  no. 5, 1085--1100.





\end{thebibliography}
\end{document}